\documentclass{article}
\usepackage{arxiv}
\usepackage{hyperref}
\usepackage[utf8]{inputenc} 
\usepackage[T1]{fontenc}    
\usepackage{url}            
\usepackage{booktabs}       
\usepackage{amsfonts}       
\usepackage{nicefrac}       
\usepackage{microtype}      
\usepackage{lipsum}		
\usepackage{graphicx}
\usepackage{doi}

\usepackage{indentfirst}
\usepackage{float}
\usepackage{tikz}
\usepackage{enumitem}           
\usepackage{amsmath,amsbsy,amsfonts,amssymb}   
\usepackage{times}
\usepackage{mathptmx} 

\title{Optimal screening strategies in the control of an infectious disease: a case of the COVID-19 in a population with age structure}

\author{ \href{https://orcid.org/0009-0007-8132-7118}{\includegraphics[scale=0.06]{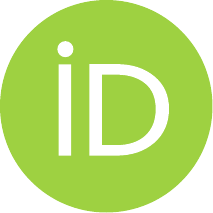}\hspace{1mm}Nelson L. Santos Junior} \\
	Department of Mathematics\\
	Federal University of Pernambuco\\
	Recife, Pernambuco, Brazil \\
	\texttt{nelson.leal@ufpe.br} \\
	\And
	\href{https://orcid.org/0000-0003-4855-4136}{\includegraphics[scale=0.06]{orcid.pdf}\hspace{1mm}João A. M. Gondim} \\
	Department of Mathematics\\
	Federal University of Pernambuco\\
	Recife, Pernambuco, Brazil \\
	\texttt{joao.gondim@ufpe.br} \\
}

\renewcommand{\shorttitle}{\textit{arXiv} Template}

\hypersetup{
pdftitle={A template for the arxiv style},
pdfsubject={q-bio.NC, q-bio.QM},
pdfauthor={David S.~Hippocampus, Elias D.~Striatum},
pdfkeywords={First keyword, Second keyword, More},
}

\begin{document}
\maketitle

\begin{abstract}
After the COVID-19 pandemic, we saw an increase in demand for epidemiological mathematical models. The goal of this work is to study the optimal control for an age-structured model as a strategy of quarantine of infected people, which is done via Pontryagin's maximum principle. Since quarantine campaigns are not just a matter of public health, also posing economic challenges, the optimal control problem does not simply minimize the number of infected individuals. Instead, it jointly minimizes this number and the economic costs associated to the campaigns, providing data that can help authorities make decisions when dealing with epidemics. The controls are the quarantine entrance parameters, which are numerically calculated for different lengths of isolation. The best strategies gives a calendar that indicates when the isolation measures can be relaxed, and the consequences of a delay in the start of the quarantine are analyzed by presenting the reduction in the number of deaths for the strategy with optimal control compared to a no-quarantine landscape.
\end{abstract}


\section{Introduction}

The mathematical models have become an important tool for investigations and predictions of infectious disease dynamics in recent decades. These works help understand how different mitigation strategies work when dealing with an epidemic \cite{mathematicalmodels}. These models gained even more interest since the last pandemic.

The COVID-19 is a infectious disease with first cases reported in Wuhan, Hubei Province in China, in December 31, 2019 \cite{whoreport1}. The World Health Organization (WHO) declared this disease a global pandemic on 11 March 2020 \cite{whoreport51}. This pandemic was one of the factors that motivated a greater demand for mathematical studies in epidemiology. From 2020 onwards, the scientists have modeled the disease, with a goal of predicting the impact of propagation of the COVID-19 outbreak \cite{wu1,ferguson,hellewell}, as well as for knowing and understanding the virus and the disease \cite{wu2}, which can influence the decisions made by governments and organizations regarding the disease.

The virus of COVID-19 quickly spread across the Americas, Europe and other continents, prompting many countries to implement lockdown measures to contain the disease. An effective early strategy to slow down the spread of the disease consisted of lockdowns, as indicated by many studies such as \cite{barkur,domenico,pepe,lau}. In March, Brazil began adopting such quarantine measures, and in March 24, for example, the partial lockdown was ordered by the state government of São Paulo \cite{govsp}.  

Hence, it was important to find out what was the best way of implementing these quarantines, which brought Optimal Control Theory into the spotlight. This theory has been applied to models of diseases such as HIV \cite{hiv1,hiv2,hiv3}, tuberculosis \cite{tub1,tub2}, influenza \cite{gripe} and general epidemics \cite{behncke,mateus}. Since 2020, some works use the optimal control applied to COVID-19 \cite{demasse,shen,obsu,zamir}. 

Gondim and Machado \cite{Gondim2} used optimal control to analyze an SEIR model with quarantine of susceptible individuals, representing a lockdown system during the COVID-19 epidemic. However, lockdowns can be seen as overly strict and may lead to substantial economic and social side effects \cite{effect1, effect2, effect3}. While lockdowns were effective measures at the onset of the outbreak, when no other options were available to prevent as many deaths as possible, as the disease progresses and becomes endemic or less aggressive, other, less costly policies should be implemented.

This paper focuses on a SEIR model with isolation of only infected people, which is closer to the proper definition of quarantine, and can be achieved by screening of symptomatic individuals through of positive tests \cite{test1,test2}. We also divide the population into age groups as in Gondim and Machado \cite{Gondim2}, this is important since the impact of disease be can different for each group. We developed the model based on the transmission dynamics and we will develop the study using data from COVID-19 in Brazil in 2020.

Our goal is to calculate numerically the optimal controls that describe the optimal time entry into quarantine of an infected individual after testing positive for the disease, minimizing not only the number of infected people but also the economic costs associated with control. This is important because even though tests are much cheaper than lockdowns or vaccines, it is still necessary to have a large enough amount of them available every day and to spread them across the country. 

With that, we hope that our results can help decisions governments must make when implementing the quarantine policies. The optimal controls gives us a calendar that indicates when can we relax the time entry into quarantine, making it possible to relax the frequency of testing in each of the age groups. Finally, we observe how the controls reduce the number of deaths in comparison to the same period without implementation of quarantine.

\section{The age-structured SEIRQ model}

Our model consists of a classical SEIR model with a quarantine of infected individuals as soon as they receive a positive test for COVID-19. The population will be divided into three age groups, according to Table \ref{groups}.

\begin{table}[H]
\caption{\label{groups} Description of the groups. }
\begin{center}
\begin{tabular}{c|c}
\hline
\textbf {Group} & \textbf{Description} \\
\hline
1 & Young people, aged 0 to 19 \\
2 & Adults, aged 20 to 59 \\
3 & Elderly, aged 60 onwards \\		     
\hline
\end{tabular}
\end{center}
\end{table}

Let $S_i(t), E_i(t), I_i(t), R_i(t)$ and $Q_i(t)$ be the number of susceptible, exposed, infected, recovered and quarantined individuals in each age group at time $t \geq 0$, respectively. We assume that the total population
\begin{equation*}
N(t)=\sum_{i=1}^3\big(S_i(t)+E_i(t)+I_i(t)+R_i(t)+Q_i(t)\big)
\end{equation*}
remains constant. That is a reasonable assumption since we are only dealing with a short time frame in comparison to the demographic time scale.

The model's equations, for $i=1,2,3$, are described in \eqref{model}.
\begin{equation}\label{model}
\left\{\begin{array}{l}
\displaystyle \frac{d S_i}{dt} = -\frac{S_i}{N}\left(\sum_{j=1}^3 \beta_{i j} I_j\right) \\
\displaystyle \frac{d E_i}{dt} = \frac{S_i}{N}\left(\sum_{j=1}^3 \beta_{i j} I_j\right)-\sigma_i E_i \\
\displaystyle \frac{d I_i}{dt} = \sigma_i E_i-\gamma_i I_i-u_i(t) I_i \\
\displaystyle \frac{d R_i}{dt} = \gamma_i I_i + \tau Q_i \\
\displaystyle \frac{d Q_i}{dt} = u_i(t) I_i - \tau Q_i
\end{array}\right.
\end{equation}

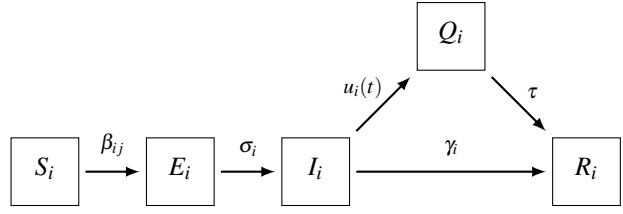
\begin{figure}[H]
\centering
\begin{tikzpicture}[scale=0.45]

\draw[black] (-9, -1) rectangle (-7,1) {}; 
\node[circle,inner sep=0.7pt,label=above:{$S_i$}] at (-8,-0.6) {};


\draw[black] (-5, -1) rectangle (-3,1) {}; 
\node[circle,inner sep=0.7pt,label=above:{$E_i$}] at (-4,-0.6) {};

\draw[thick][-latex] (-6.8,0) -- (-5.2,0);
\node[circle,inner sep=0.7pt,label=above:{ \scalebox{0.8}{$\beta_{ij}$}}] at (-6, 0.1) {};


\draw[black] (-1,-1) rectangle (1,1); 
\node[circle,inner sep=0.7pt,label=above:{$I_i$}] at (0,-0.6) {};

\draw[thick][-latex] (-2.8,0) -- (-1.2,0);
\node[circle,inner sep=0.7pt,label=above:{ \scalebox{0.8}{$\sigma_i$}}] at (-2,0.1) {};

 
\draw[black] (3,3) rectangle (5,5);
\node[circle,inner sep=0.7pt,label=above:{$Q_i$}] at (4,3.4) {};

\draw[thick][-latex] (1.2,1.2) -- (2.8,2.8);
\node[circle,inner sep=0.7pt,label=above:{ \scalebox{0.8}{$u_i(t)$}}] at (1.4,1.8) {};


\draw[thick][-latex] (5.2,2.8) -- (6.8,1.2);
\node[circle,inner sep=0.7pt,label=above:{ \scalebox{0.8}{$\tau$}}] at (6.4,1.8) {};
 
\draw[black] (7,-1) rectangle (9,1) {}; 
\node[circle,inner sep=0.7pt,label=above:{$R_i$}] at (8,-0.6) {};

\draw[thick][-latex] (1.2,0) -- (6.8,0);
\node[circle,inner sep=0.7pt,label=above:{ \scalebox{0.8}{$\gamma_i$}}] at (4,0.1) {};

\end{tikzpicture}
\vspace*{0.3cm}\caption{\label{diagram} Compartment diagram of the model.}
\end{figure}

Figure \ref{diagram} represents the dynamics between the compartments of our model.

All parameters are non-negative. $\beta_{ij}$ is the transmission coefficient from the age group $i$ to the age group $j$. As in \cite{castilho}, it will be assumed that $\beta_{ij} = \beta_{ji}$ for all $i$ and $j$, $\sigma_i$ and $\gamma_i$ are the latency and recovery periods, respectively, for age group $i$, $\tau$ is the exit rate from quarantine. Parameter values are taken from \cite{castilho}, where data fitting was performed using an adaptation of a least squares algorithm from \cite{martcheva}, and are listed in Table \ref{parameters}.

\begin{table}[H]
\caption{\label{parameters} Parameters values. }
\begin{center}
\begin{tabular}{c|c|c|c}
\hline
\textbf {Parameter} & \textbf{Value} & \textbf {Parameter} & \textbf{Value}\\
\hline
$\beta_{11}$ & $1.76168$ & $\sigma_{1}$ & $0.27300$ \\
$\beta_{12}$ & $0.36475$ & $\sigma_{2}$ & $0.58232$ \\
$\beta_{13}$ & $1.32468$ & $\sigma_{3}$ & $0.69339$ \\	     
$\beta_{22}$ & $0.63802$ & $\gamma_{1}$ & $0.06862$ \\
$\beta_{23}$ & $0.35958$ & $\gamma_{2}$ & $0.03317$ \\
$\beta_{33}$ & $0.57347$ & $\gamma_{3}$ & $0.35577$ \\
\hline
\end{tabular}
\end{center}
\end{table}

In the WHO's recommendation of May 27, 2020 \cite{oms}, the criteria for discharge from isolation of a COVID-19 patient is 10 days after symptom onset plus at least 3 additional days without symptoms for symptomatic patients and 10 days after a positive SARS-CoV-2 test for asymptomatic cases. We consider that an infected individual stays in quarantine, then, for 13 days, therefore the quarantine exit parameter is taken as $\tau=\dfrac{1}{13}$. 

Entry into quarantine will occur after a positive test for the disease. The $u_i(t)$ denotes the controls, which represent the fraction of infected individuals in each age group that are quarantined per unit of time at time $t$, after the positive result. The controls satisfy, therefore,
\begin{equation}\label{limitcontrol}
0 \leq u_i(t) \leq 1, \quad i=1,2,3, 
\end{equation}
since no will exist maximum time for  to entry quarantine and the minimum time is 1 day.

The total population of group $i$ will be defined as
\begin{equation*}    
N_i(t) = S_i(t) + E_i(t) + I_i(t) + R_i(t) + Q_i(t).
\end{equation*}

By adding the equations in system \eqref{model}, it is clear that $N_i(t)$ is also constant for $i=1,2,3$. Moreover, $R_i(t)$ only appears in the other equations as a part of $N_i(t)$, so we substitute the equations for $R_i'(t)$ by $N_i'(t)$. Hence, we are able to rewrite the system of equations as in \eqref{statesystem}.

\begin{equation}\label{statesystem}
\left\{\begin{array}{l}
\displaystyle \frac{d S_i}{dt} = -\frac{S_i}{N}\left(\sum_{j=1}^3 \beta_{i j} I_j\right) \\
\displaystyle \frac{d E_i}{dt} = \frac{S_i}{N}\left(\sum_{j=1}^3 \beta_{i j} I_j\right)-\sigma_i E_i \\
\displaystyle \frac{d I_i}{dt} = \sigma_i E_i-\gamma_i I_i-u_i(t) I_i \\
\displaystyle \frac{d Q_i}{dt} = u_i(t) I_i-\tau Q_i \\
\displaystyle \frac{d N_i}{dt} = 0 
\end{array}\right.
\end{equation}


We suppose, as a simplification, that the number of recoveries is given by the difference between the total number of cases and the number of deaths. Even though this might not reflect the reality, since there could be some patients who we considered as recovered but still carry the disease, we proceed this way nonetheless due to the scarcity of available information on recoveries. 


The initial time for our simulations will be set as May 8, 2020, and we use data from Brazil to set the initial conditions. According to \cite{Worldmeters,Gondim2}, there were 65,124 recovered and 76,603 active cases in the country. The mean incubation period of the disease has been estimated to be around $5$ days \cite{lauer}. Also according to the same references, there were 97,575 active COVID-19 cases in Brazil by May 13, 2020, so we set the initial number of exposed individuals as the increment of 20,972 active cases from these $5$ days.  We assume that there are no quarantined individuals at the start of the simulations and, as the numbers of exposures, infections and recoveries are rather small relative to the total population, we assume that the number of susceptible individuals for each age group is the same as the total population of that respective age group. A summary of the initial conditions, rounded to the nearest integers, is shown in Table \ref{initialconditions}.


Using the data available in a 2020 report from Spain's Ministry of Health \cite{espanha}, we calculate how infections and deaths are distributed among the age groups, and the percentages are available in Table \ref{rateespanha}. Finally, we fix the total population at 200 million, of which 40\% are young individuals, 50\% are adults and 10\% are elderly. 

\begin{table}[H]
\caption{\label{rateespanha} Number of cases, recoveries and deaths for each group by \cite{espanha}.}
\begin{center}
\begin{tabular}{c|c|c|c|c|c}
\hline
\textbf {Group} & \textbf{Cases} &  $\%$ \textbf{of cases} & \textbf{Recoveries} & $\%$ \textbf{of recoveries} & \textbf{Deaths}\\
\hline
1 & 2448 & $1.03\%$ &  2441 & $1.12\%$ & 7 \\
2 & 113059 & $47.62\%$ & 112168 & $51.31\%$ & 981 \\
3 & 121928 & $51.35\%$ & 103980 & $47.57\%$ & 17948 \\
\hline
Total& 237435 & $100\%$ & 218589 & $100\%$ & 18846 \\
\hline
\end{tabular}
\end{center}
\end{table}

\begin{table}[H]
\caption{\label{initialconditions} Initial conditions }
\begin{center}
\begin{tabular}{c|c|c|c}
\hline
\textbf {Class} & \textbf{$i=1$} & \textbf {$i=2$} & \textbf{$i=3$}\\
\hline
Susceptible & 80000000 & 10000000 & 20000000 \\
Exposed & 216 & 9987 & 10769 \\
Infected & 789 & 36478 & 39335 \\	     
Recovered & 729 & 33415 & 30979 \\
Quarantined & 0 & 0 & 0 \\
\hline
\end{tabular}
\end{center}
\end{table}

\section{The optimization problem}

Following system (\ref{statesystem}) , we minimize the functional 
\begin{equation}\label{functional}
J = \int_0^T \sum_{i=1}^3 \big( I_i(t) + B_i u^2_i(t) \big)~dt. 
\end{equation}

In \eqref{functional}, $T$ is the quarantine duration, while the parameters $B_i$ are the costs of the control. It will be assumed that $B_i>0$ for $i=1,2,3$, as well as that 
\begin{equation}
B = B_1 + B_2 + B_3,   
\end{equation}
where $B \in \mathbb{R}$ is the total control cost.

Sufficient conditions for the existence of optimal controls are available in standard results from optimal control theory. In this case, Theorem 2.1 in Joshi et al. \cite{Joshi} can be used to show that the optimal control exists. 

Now, Pontryagin's maximum principle \cite{pontryagin,lenhart} states that the optimal controls are solutions of the Hamiltonian system
\begin{equation}\label{hamiltonian}
\begin{aligned}
\displaystyle H=&\sum_{i=1}^3\left[B_i u_i^2+u_i I_i\left(\lambda_Q^i-\lambda_I^i\right)+\frac{S_i}{N}\left(\sum_{j=1}^3 \beta_{ij} I_j\right)\left(\lambda_E^i-\lambda_S^i\right)\right] \\
&+ \sum_{i=1}^3\left[\sigma_i E_i\left(\lambda_I^i-\lambda_E^i\right)+I_i\left(1-\gamma_i \lambda_I^i\right)-\tau Q_i \lambda_Q^i \right],
\end{aligned}
\end{equation}
involving the state variables $S_i, E_i, I_i, R_i, Q_i$ and the adjoint variables $\lambda_S^i$, $\lambda_E^i$, $\lambda_I^i$, $\lambda_Q^i$, $\lambda_N^i$, which satisfy the adjoint system 
\begin{equation}\label{adjointsystem}
\left\{\begin{array}{l}
\displaystyle \frac{d\lambda_S^i}{dt}=\frac{1}{N}\left(\sum_{j=1}^3 \beta_{ij} I_j\right)\left(\lambda_S^i-\lambda_E^i\right) \\
\displaystyle \frac{d\lambda_E^i}{dt}=\sigma_i\left(\lambda_E^i-\lambda_I^i\right) \\
\displaystyle \frac{d\lambda_I^i}{dt}=u_i\left(\lambda_I^i-\lambda_Q^i\right)+\gamma_i \lambda_I^i-1+\frac{1}{N}\left(\sum_{j=1}^3 \beta_{ij} S_j\left(\lambda_S^j-\lambda_E^j\right)\right) \\
\displaystyle \frac{d\lambda_Q^i}{dt}=\tau \lambda_Q^i \\
\displaystyle \frac{d\lambda_N^i}{dt}=\frac{1}{N^2}\left( \sum_{j,k=1}^3 \beta_{k j} S_k I_j\left(\lambda_E^k-\lambda_S^k\right)\right)
\end{array}\right.   
\end{equation}

The adjoint variables also must satisfy the transversality conditions
\begin{equation}\label{transversalitycondition}
\lambda_S^i(T)=\lambda_E^i(T)=\lambda_I^i(T)=\lambda_Q^i(T)=\lambda_N^i(T)=0, \quad i=1,2,3. 
\end{equation}

Finally, the optimal controls are critical point of Hamiltonian function, which gives the optimality condition
\begin{equation*}
\frac{\partial H}{\partial u_i}\Bigg|_{u_i=u_i^*}=0, 
\end{equation*}
resulting in
\begin{equation}
\displaystyle u_i^*=\frac{I_i\left(\lambda_I^i-\lambda_Q^i\right)}{2 B_i}. 
\end{equation}

Since we are considering bounded controls, by (\ref{limitcontrol}), the $u_i^*$ are calculated \cite{lenhart} using
\begin{equation}\label{optimalycondition}
\displaystyle u_i^*=\min \left\{1, \max \left\{0, \frac{I_i\left(\lambda_I^i-\lambda_Q^i\right)}{2 B_i}\right\}\right\}, 
\end{equation}

Uniqueness of the optimal controls (at least for small enough $T$) also follow from standard results, such as Theorem 2.3 in Joshi et al. \cite{Joshi} .  

The numerical solutions of systems (\ref{statesystem}) and (\ref{adjointsystem}) can be found using the forward-backward sweep method \cite{lenhart}.  We have a initial condition in $t=0$ for system (\ref{statesystem}) and a initial condition in $t=T$ for system (\ref{adjointsystem}), the algorithm of sweep uses these conditions and works with the following steps:

\textbf{Step 1}: Starts with an initial guess of the controls $u_1$, $u_2$ and $u_3$, solves (\ref{statesystem}) using the initial condition for state variable and the values for $u_i$ forward in time;

\textbf{Step 2}: Uses initial condition (\ref{transversalitycondition}) and the results provided in step 1 for the state variables and for the optimal control, to solve (\ref{adjointsystem}) backward in time, the new controls are defined following (\ref{optimalycondition});

\textbf{Step 3}: This process continues until a convergence is obtained.

\section{Discussion}

Infected individuals are identified by screening carried out by tests. The purchase, storage and distribution of tests, in addition to the organization of medical teams and advertising campaigns for testing, are of significant importance in these containment strategies. In this regard, quarantine is not just a matter of public health; it also involves economic considerations. Therefore, each group $i$ is associated with a control cost $B_i$.

According to \cite{castilho}, the sensitivity coefficient $R_0$ exhibits a greater decay for group 1 with an increase in the recovery coefficient $\gamma_i$; this is why class 1 is the most sensitive to screening measures, indicating that youngsters should be preferentially screened. This implies that group 1 requires a greater investment in testing, resulting in a higher cost. Furthermore, as $\gamma_i$ increases,$R_0$ shows similar behavior for groups 2 and 3 \cite{castilho}. Although the number of tests is higher for group 2 due to its larger population, group 3 requires greater investment in a publicity campaign because it has limited access to information, and we assume these groups have the same control cost. Therefore, we allocate a total control cost distributed as 40\%, 30\%, and 30\% for groups 1, 2, and 3, respectively.

\begin{figure}[H]
\begin{center}
\includegraphics[scale=0.9, trim=20mm 115mm 20mm 115mm]{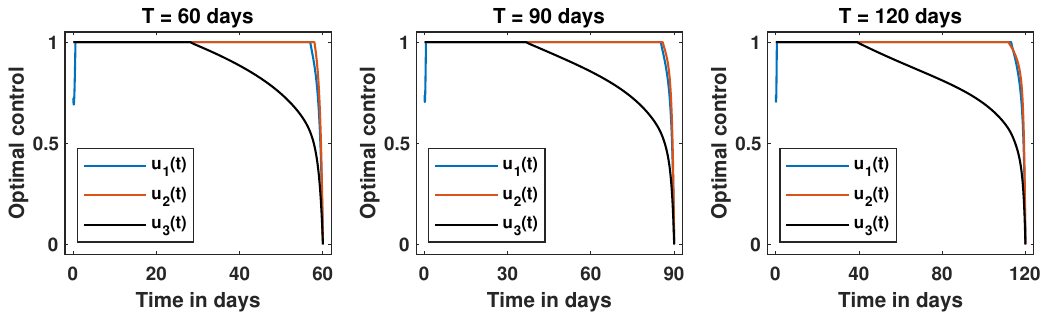}
\caption{\label{controls}The optimal controls for different quarantine lengths.}
\end{center}
\end{figure}

Figure \ref{controls} plots the graphs of the optimal controls $u_1(t)$, $u_2(t)$ and $u_3(t)$ for the cost distribution. The simulation is done to 60, 90 and 120 days. 

An interesting feature of the plots of optimal control is that they provide an``optimal calendar'' of when the quarantine should start to be relaxed. This calendar is shown in Table \ref{calendar1}. 

\begin{table}[H]
\caption{\label{calendar1} Description of the optimal calendar from optimal control.}
\begin{center}
\begin{tabular}{c|c|c|c} 
\hline
\textbf {Group} & \textbf{$T=60$} & \textbf{$T=90$} & \textbf{$T=120$}  \\
\hline
1 & $57$ days & $85$ days & $113$ days \\
\hline
2 & $58$ days & $86$ days & $112$ days \\
\hline
3 & $28$ days & $37$ days & $39$ days \\		     
\hline
\end{tabular}
\end{center}
\end{table}

According to \cite{Gondim2}, for a lockdown, the relaxation of isolation happens as follows: adults, who incur higher costs, are freed first, and the elderly, who have the lowest costs, are freed last. In other words, the relaxation follows a decreasing order of population distribution. In our model, which considers only the quarantine of infected individuals, as shown in Table \ref{calendar1}, the order of relaxation is the inverse: the elderly, who have the smallest population, are liberated first, and adults with the largest population are liberated last. The cost distribution we used is different; however, even if we adopted the cost distribution from \cite{Gondim2}, there would be no change in the behavior of the optimal control regarding the order of relaxation between the groups. We believe this difference arises because the screening process is distinct in each case.

\begin{figure}[H]
\begin{center}
\includegraphics[scale=0.9, trim=20mm 115mm 20mm 115mm]{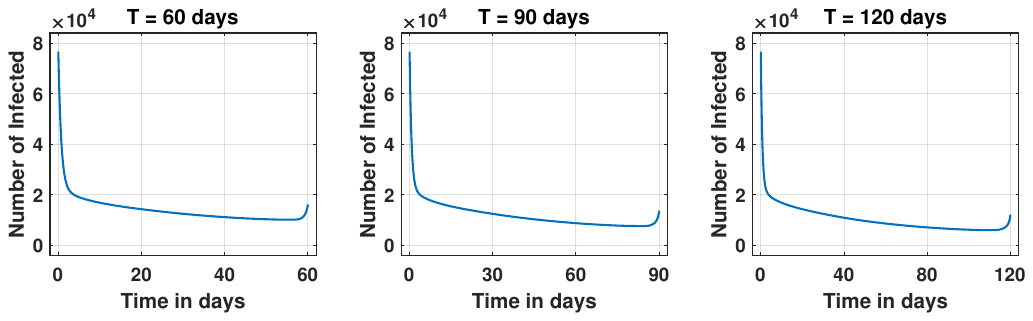}
\caption{\label{infectedcontrol} Curves of infections for the optimal control in different quarantine lengths.}
\end{center}
\end{figure}

Figure \ref{infectedcontrol} plots the graphs of the infected curves for an optimal control quarantine. In all cases of quarantine duration, the number of infected reaches a minimum just before the end of the simulation, and by the time the quarantine ends, the number of cases is increasing but remains much lower than the initial count.

\begin{figure}[H]
\begin{center}
\includegraphics[scale=0.9, trim=20mm 110mm 20mm 110mm]{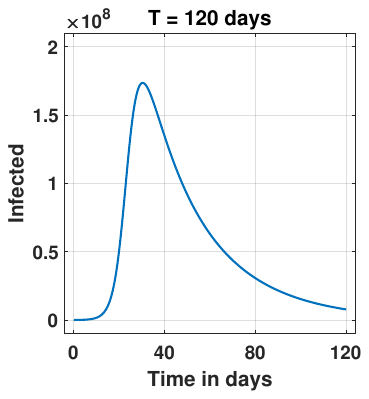}
\caption{\label{infectednocontrol} Curve of infections without quarantine for 120 days.}
\end{center}
\end{figure}

The importance of quarantine is reflected in the number of infected, which is shown in Figure \ref{infectednocontrol}, where we did not consider a control of quarantine and we have a large cumulative number of infected. The curves of infected without control of quarantine for the 60 and 90 days are not presented because they behave like restrictions of the curve for 120 days.

We can analyze what happens if the government takes a long time to implement the quarantine. As of 13 May 2020, the number of active cases in Brazil doubles after 10 days and quadruples after 20 days \cite{Gondim2}. We do this by considering the initial conditions of exposed, infected, and recovered individuals to be twice and four times as large as their original values.

\begin{figure}[H]
\begin{center}
\includegraphics[scale=0.9, trim=20mm 115mm 20mm 115mm]{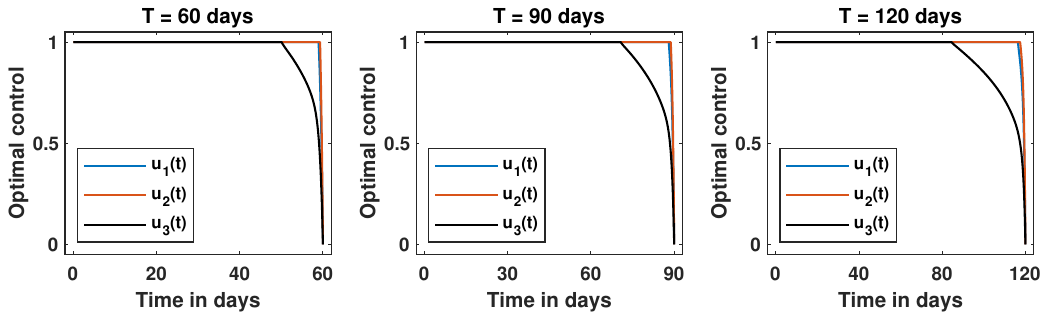}
\caption{\label{optimalcontrol2} Optimal control with 10-day delay at the start of the quarantine.}
\end{center}
\end{figure}

\begin{figure}[H]
\begin{center}
\includegraphics[scale=0.9, trim=20mm 115mm 20mm 115mm]{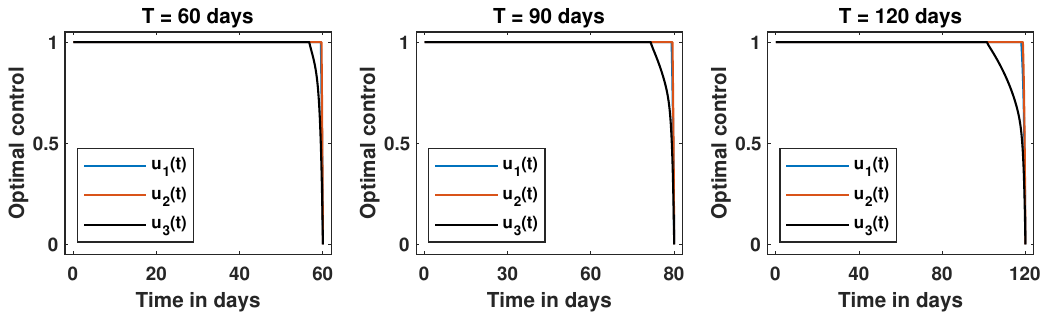}
\caption{\label{optimalcontrol3} Optimal control with 20-day delay at the start of the quarantine.}
\end{center}
\end{figure}

The controls shown in Figures \ref{optimalcontrol2} and \ref{optimalcontrol3} provide different relaxation calendars, which are described in Tables \ref{calendar2} and \ref{calendar3}.

\begin{table}[H]
\caption{\label{calendar2} Optimal calendar with delay of 10 days in start of the quarantine. }
\begin{center}
\begin{tabular}{c|c|c|c} 
\hline
\textbf {Group} & \textbf{$T=60$} & \textbf{$T=90$} & \textbf{$T=120$}  \\
\hline
1 & $59$ days & $88$ days & $116$ days \\
\hline
2 & $59$ days & $89$ days & $117$ days \\
\hline
3 & $50$ days & $71$ days & $84$ days \\		     
\hline
\end{tabular}
\end{center}
\end{table}

\begin{table}[H]
\caption{\label{calendar3} Optimal calendar with delay of 20 days in start of the quarantine. }
\begin{center}
\begin{tabular}{c|c|c|c} 
\hline
\textbf {Group} & \textbf{$T=60$} & \textbf{$T=90$} & \textbf{$T=120$}  \\
\hline
1 & $59$ days & $89$ days & $118$ days \\
\hline
2 & $60$ days & $89$ days & $119$ days \\
\hline
3 & $57$ days & $81$ days & $101$ days \\		     
\hline
\end{tabular}
\end{center}
\end{table}

Comparing the calendars of Tables \ref{calendar1}, \ref{calendar2}, and \ref{calendar3}, the delay of 10 and 20 days in the start of isolation of the infected results in a small extension of the quarantine time. This extension is sufficient for there to be practically no relaxation for groups 1 and 2. For group 3, the quarantine time increases by 0.78 and 1.03 times for a quarantine of 60 days, by 0.91 and 1.18 times for a quarantine of 90 days, and by 1.15 and 1.58 times for a quarantine of 120 days.

We will analyze how a quarantine with optimal control reduce the number of deaths. In our model, deaths are calculated as a fraction of the recovered individuals. From Table \ref{rateespanha}, we have the fatality rates $\mu_1$, $\mu_2$, and $\mu_3$ for groups 1, 2, and 3, respectively. The rates are
\begin{equation*}
\mu_1 = 0.003 \quad \mu_2 = 0.008 \quad \mu_3 = 0.147.
\end{equation*}

\begin{figure}[H]
\begin{center}
\includegraphics[scale=0.9, trim=20mm 115mm 20mm 115mm]{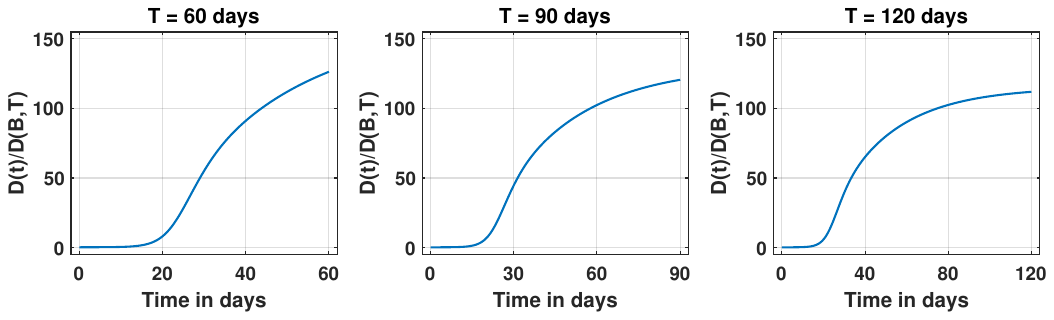}
\caption{\label{deaths} Plots of $D(t)/D(B,T)$ for different quarantine lengths.}
\end{center}
\end{figure}

Let's denote $D(B,T)$ as the cumulative number of deaths in an optimal control campaign with cost $B$ and final time $T$, and $D(t)$ as the number of deaths for the model with no quarantine at time $t$. In Figure \ref{deaths}, we plot graphs of $D(t)$ divided by $D(B,T)$ for $T = 60, 90 \text{ and } 120$ days.

Due to the uncertain nature of the parameters and to the large number of unreported cases, we do not show the crude numbers of $D(B,T)$ and $D(t)$. However, at the end of the quarantine, optimal controls reduce the number of deaths by 126.4, 120.5, and 111.9 times for 60, 90, and 120 days, respectively.

All simulations are performed using the Matlab 2018 techniques \cite{matlab}.

\section{Conclusions}

In this paper, we consider an age-structured SEIRQ model, where the entry parameters of infected individuals in quarantine are treated as controls of the system. We determined the optimal controls using Pontryagin’s maximum principle and employed the forward-backward sweep method to calculate the optimal controls numerically.

We analyze how control influences results, including relaxation calendars, the number of infections, and the number of deaths at the end of quarantine. The division of costs and the results for optimal controls provide data that can help authorities make difficult decisions, such as when to start or relax quarantine measures, while minimizing economic impact for each age group.

The calendar for relaxation of isolation measures in the three age groups, based on quarantine lengths of 60, 90, and 120 days, shows a reasonable relaxation period for the elderly and a very short period for young people and adults. Regarding the order of relaxation, the elderly are liberated first and adults are liberated last, an inverse order to the lockdown case in \cite{Gondim2}.

We show that with a delay of 10 and 20 days in the start of quarantine, the relaxation period for groups 1 and 2, which was already very limited, becomes practically non-existent. For group 3, the quarantine duration increases by 0.78 and 1.03 times for a 60-day quarantine, by 0.91 and 1.18 times for a 90-day quarantine, and by 1.15 and 1.58 times for a 120-day quarantine.

In all investigated cases, the number of infected individuals reached a minimum just before the end of the simulation, and the cases began to rise again at the end of isolation, although they reached a much lower number than initially. This shows that the quarantine was effective but could have had a longer duration.

Finally, compared to the same period of time without quarantine, the optimal controls produced a reduction in the number of deaths in 126.4, 120.5 and 111.9 times, for 60, 90 and 120 days, respectively.

In our model we used data from Brazil as initial conditions and in the parameter fitting. Brazil is a very large country, with many cities at different stages of the pandemic. This means that studies such as this one should be made locally to best suit the characteristics of each city or state.

This work was developed with dynamics and data related to COVID-19. The results can be used for a potential new outbreak of this disease or can serve as a structured basis for studies on other diseases. For this, we must adapt the dynamics and perform a good estimate of the parameters.



\end{document}